\documentclass[12pt]{article}

\usepackage[english,francais]{babel}
\usepackage{latexsym}
\usepackage{amssymb}
\usepackage{amsfonts}
\usepackage{pstricks}
\usepackage{amsmath}

\newtheorem{thm}{Theorem}
\newtheorem{lem}{Lemma}
\newtheorem{cor}{Corollary}
\newtheorem{rem}{Remark}
\newtheorem{prop}{Proposition}

\newtheorem{Fact}{Fact}

\newcommand{\bpr}{\noindent{\bf Proof\/. }}
\newcommand{\epr}{\hspace*{\fill}$\square$\medskip}

\author{Houmem Belkhechine\\ Universit\'e de Carthage\\
Institut Pr\'eparatoire aux \'Etudes d'Ing\'enieurs de Bizerte, Tunisie\\
{\tt houmem@gmail.com} \and
Imed Boudabbous\\ Universit\'e de Sfax\\
Institut Pr\'eparatoire aux \'Etudes d'Ing\'enieurs de Sfax, Tunisie\\
{\tt imed.boudabbous@gmail.com}\and
Kaouthar Hzami\\ Universit\'e de Gab\`es\\
Institut sup\'erieur d'Informatique et de Multim\'edia de Gab\`es, Tunisie\\
{\tt hzamikawthar@gmail.com}}

 \title{The indecomposable tournaments $T$ with $\mid\! W_{5}(T) \!\mid = \mid\! T \!\mid -2$}

\date{}

\begin{document}

\maketitle
\selectlanguage{english}
\begin{abstract}
We consider a tournament $T=(V, A)$. For $X\subseteq V$, the subtournament of $T$ induced by $X$ is $T[X] = (X, A \cap (X \times X))$. An interval of $T$ is a subset $X$ of $V$ such that for $a, b\in X$ and $ x\in V\setminus X$, $(a,x)\in A$ if and only if $(b,x)\in A$.
The trivial intervals of $T$ are $\emptyset$, $\{x\}(x\in V)$ and $V$. A tournament is indecomposable if all its intervals are trivial. For $n\geq 2$, $W_{2n+1}$ denotes the unique indecomposable tournament defined on $\{0,\dots,2n\}$ such that $W_{2n+1}[\{0,\dots,2n-1\}]$ is the usual total order. Given an indecomposable tournament $T$, $W_{5}(T)$ denotes the set of $v\in V$ such that there is $W\subseteq V$ satisfying $v\in W$ and $T[W]$ is isomorphic to $W_{5}$. Latka \cite{BJL} characterized the indecomposable tournaments $T$ such that $W_{5}(T)=\emptyset$. The authors \cite{HIK} proved that if $W_{5}(T)\neq \emptyset$, then $\mid\! W_{5}(T) \!\mid \geq \mid\! V \!\mid -2$. In this article, we characterize the indecomposable tournaments $T$ such that $\mid\! W_{5}(T) \!\mid = \mid\! V \!\mid -2$.

\end{abstract}
\section{Introduction}
\subsection{Preliminaries}
A {\it tournament} $T = (V(T), A(T))$ (or $(V, A)$) consists of a finite set $V$ of {\it vertices} together with a set $A$ of
ordered pairs of distinct vertices, called {\it arcs}, such that for all $x \neq y \in V$, $(x,y) \in A$ if and only if
$(y,x) \not\in A$. The {\it cardinality} of $T$, denoted by $\mid\!T\!\mid$, is that of $V(T)$.
Given a tournament $T=(V,A)$, with each subset $X$ of $V$ is associated the {\it subtournament} $T[X] = (X, A \cap (X \times X))$
of $T$ induced by $X$. For $X\subseteq V$ (resp. $x\in V$), the subtournament $T[V \setminus X]$
(resp. $T[V \setminus \{x\}$]) is denoted by $T-X$ (resp. $T-x$).
Two tournaments $T = (V, A)$ and $T' = (V', A')$ are {\it isomorphic}, which is denoted by $T \simeq T'$, if there exists an {\it isomorphism} from $T$ onto $T'$, i.e., a bijection $f$ from $V$ onto $V'$ such that for all $x$, $y \in V$, $(x, y) \in A$ if and only if $(f(x), f(y)) \in A'$.
We say that a tournament $T'$ {\it embeds} into $T$ if $T'$ is isomorphic to a subtournament of $T$.
Otherwise, we say that $T$ {\it omits} $T'$. The tournament $T$ is said to be
{\it transitive} if it omits the tournament $C_{3} = (\{0,1,2\}, \{(0,1), (1,2), (2,0)\})$.
For a finite subset $V$ of $\mathbb{N}$, we denote by $\overrightarrow{V}$ the usual {\it total order} defined on $V$, i.e., the transitive tournament $(V, \{(i,j): i<j\})$.

Some notations are needed. Let $T=(V,A)$ be a tournament.
For two vertices $x \neq y \in V$, the notation $x \longrightarrow y$ signifies that $(x, y) \in A$.
Similarly, given $x\in V$ and $Y\subseteq V$, the notation $x\longrightarrow Y$ (resp. $Y \longrightarrow x$)
means that $x\longrightarrow y$ (resp. $y \longrightarrow x$) for all $y \in Y$. Given $x\in V$, we set $N_{T}^{+}(x)=\{y \in V$: $x \longrightarrow y \}$.
For all $n\in \mathbb{N}\setminus\{0\}$, the set $\{0,\dots, n-1\}$ is denoted by $\mathbb{N}_{n}$.

Let $T=(V,A)$ be a tournament. A subset $I$ of $V$ is an {\it interval} \cite{F1,I,ST} (or a {\it clan} \cite{E})
of $T$ provided that for all $x \in V \setminus I$, $x \longrightarrow I$ or $I\longrightarrow x$.
For example, $\emptyset $, $\{x\}$, where $x \in V$, and $V$ are intervals of $T$, called {\it trivial} intervals. A tournament
is {\it indecomposable} \cite{I,ST} (or {\it primitive} \cite{E})
if all its intervals are trivial, otherwise it is {\it decomposable}.
Notice that a tournament $T=(V,A)$ and its {\it dual} $T^{\star}=(V,\{(x,y): (y,x)\in A\})$ share the same intervals. Hence, $T$ is indecomposable if and
only if $T^{\star}$~is.

\begin{sloppypar}
For $n\geq2$, we introduce the tournament $W_{2n+1}$ defined on $\mathbb{N}_{2n+1}$ as follows: $W_{2n+1}[\mathbb{N}_{2n}] = \overrightarrow{\mathbb{N}_{2n}}$ and $N_{W_{2n+1}}^{+}(2n) = \{2i: i \in \mathbb{N}_{n}\}$ (see Figure~1).
In 2003, B.J. Latka  \cite{BJL} characterized the indecomposable tournaments omitting the tournament $W_{5}$. In 2012, the authors were interested in the set $W_{5}(T)$ of the vertices $x$ of an indecomposable tournament $T=(V,A)$ for which there exists a subset $X$ of $V$ such that $x \in X$ and $T[X]\simeq W_{5}$.
They obtained the following.
\end{sloppypar}
\begin{thm} [\cite{HIK}] \label{HIK}
Let $T$ be an indecomposable tournament into which $W_{5}$ embeds. Then, $\mid\! W_{5}(T) \!\mid \geq \mid\!T\!\mid-2$. If, in addition, $\mid\!T\!\mid$ is even, then $\mid\! W_{5}(T) \!\mid  \geq \mid\!T\!\mid-1$.
\end{thm}

In this paper, we characterize the class $\mathcal{T}$ of the indecomposable tournaments $T$ on at least 3 vertices such that $\mid\! W_{5}(T) \!\mid =\mid\!T\!\mid-2$. This answers [1, Problem 4.4]
\begin{center}
\psset{unit=0.85cm}
\begin{pspicture}(11,4)

\psdots (0,0) (1.5,0) (4.5,0) (6,0) (9,0) (10.5,0) (5.25,2.5)

\put (5.1,2.6){\scriptsize $2n$}

\put (-0.1, -0.4){\scriptsize ${0}$} \put (1.4, -0.4){\scriptsize $1$}

\put (4.25, -0.34) {\scriptsize $2i$} \put (5.2, -0.34){\scriptsize $2i+1$}

\put (8.2,-0.4){\scriptsize $2n-2$} \put (10, -0.4){\scriptsize $2n-1$}

\put(2.8,0) {$\ldots$} \psline{->}(1.5,0) (2.3,0)

\psline{->}(3.7,0) (4.5,0) \psline{->}(6,0) (6.8,0)

\psline{->}(8.2,0) (9,0) \put (7.3,0){$\ldots$}

\psline{->}(0,0)(0.85,0) \psline (0.75,0) (1.5,0)

\psline{->}(4.5,0)(5.35,0) \psline (5.3,0) (6,0)

\psline{->}(9,0)(9.85,0) \psline(9.8,0)(10.5,0)
\psline(0,0)(0.5,-0.4)\psline{->}(0.5,-0.4) (1,-0.4)
\psline(1.5,0)(2,-0.4)\psline{->}(2,-0.4) (2.5,-0.4)
\psline(4.5,0)(5,-0.4)\psline{->}(5,-0.4) (5.5,-0.4)
\psline(6,0)(6.5,-0.4)\psline{->}(6.5,-0.4) (7,-0.4)

\psline{->}(5.25,2.5) (2.625, 1.25) \psline(2.625,1.25) (0,0)

\psline{->}(5.25,2.5)(4.875,1.25) \psline (4.875,1.25)(4.5,0)

\psline{->}(1.5,0)(3.375,1.25) \psline (3.375,1.25) (5.25,2.5)

\psline{->}(6,0)(5.625,1.25) \psline (5.625,1.25) (5.25, 2.5)

 \psline{->}(10.5,0)(7.875,1.25)  \psline (7.875,1.25) (5.25,2.5)

\psline{->}(5.25,2.5) (7.125,1.25) \psline (7.125,1.25) (9,0)

\put (4, -1.2) {{\bf Fig.1.} $W_{2n+1}$}
\end{pspicture}

\end{center}

\vspace{1cm}

\subsection{Partially critical tournaments and the class $\mathcal{T}$}
Our characterization of the tournaments of the class $\mathcal{T}$ requires the study of their partial criticality structure, notion introduced as a weakening of the notion of critical indecomposability defined in Section \ref{sec:2}. These notions are defined in terms of critical vertices. A vertex $x$ of an indecomposable tournament $T$ is {\it critical} \cite{ST} if $T-x$ is decomposable. The set of non-critical vertices of an indecomposable tournament $T$ was introduced in \cite{MYS}. It is called the {\it support} of $T$ and is denoted by $\sigma(T)$. Let $T$ be an indecomposable tournament and let $X$ be a subset of $V(T)$ such that $\mid\! X  \!\mid \geq 3$, we say that $T$ is {\it partially critical according to $T[X]$} (or {\it $T[X]$-critical})  \cite{BDI} if $T[X]$ is indecomposable and if $\sigma(T)\subseteq X$. We will see that: {\it for $T\in \mathcal{T}$, $ V(T)\setminus W_{5}(T)=\sigma(T)$}. Partially critical tournaments are characterized by M.Y. Sayar in \cite{MYS}. In order to recall this characterization, we first introduce the tools used to this end.
Given a tournament $T=(V,A)$, with each subset $X$ of $V$, such that $\mid\!X\!\mid\geq3$ and
$T[X]$ is indecomposable, are associated the following subsets of $V\setminus X$.
\begin{itemize}
     \item $\langle X\rangle=\{x \in V\setminus X: x \longrightarrow X$\  or \ $X \longrightarrow x\}$.
    \item For all $u\in X$, $X(u)=\{x \in V\setminus  X: \{u, x\}$ is an interval of $T[X\cup \{x\}] \}$.
    \item Ext$(X)=\{x \in V\setminus  X: T[X\cup \{x\}]$ is indecomposable$\}$.
\end{itemize}
The family $\{X(u): u\in X\}\cup\{$Ext$(X), \langle X\rangle\}$ is denoted by $p^{T}_{X}$.
\begin{lem} [\cite{E}] \label{ER}
Let $T=(V,A)$ be a tournament and let $X$ be a subset of $V$ such that $\mid\!X\!\mid\geq3$ and $T(X)$
is indecomposable. The nonempty elements of $p^{T}_{X}$ constitute a partition of $V\setminus X$
and satisfy the following assertions.
\begin{itemize}
    \item For $u\in X$, $x\in X(u)$ and $y\in V\setminus (X\cup X(u))$, if $T[X\cup \{x, y\}]$ is decomposable,
    then $\{u, x\}$ is an interval of $T[X\cup \{x, y\}]$.
    \item For $x\in \langle X\rangle$ and $y\in V\setminus (X\cup \langle X\rangle)$, if $T[X\cup \{x, y\}]$ is decomposable,
    then $X\cup \{y\}$ is an interval of $T[X\cup \{x, y\}]$.
    \item For $x\neq y\in Ext(X)$, if $T[X\cup \{x, y\}]$ is decomposable,
    then $\{x, y\}$ is an interval of $T[X\cup \{x, y\}]$.
\end{itemize}
\end{lem}
Furthermore, $\langle X\rangle$ is divided into $X^{-} = \{x \in \langle X\rangle: x \longrightarrow X\}$ and $X^{+} = \{x \in \langle X\rangle: X \longrightarrow x\}$. Similarly, for all $u \in X$, $X(u)$ is divided into $X^{-}(u) = \{x \in X(u): x \longrightarrow u\}$ and $X^{+}(u) = \{x \in X(u): u \longrightarrow x\}$. We then introduce the family $q^{T}_{X}= \{$Ext$(X), X^{-}, X^{+}\} \cup \{X^{-}(u): u \in X\} \cup \{X^{+}(u): u\in X\}$.

A {\it graph} $G = (V(G),E(G))$ (or $(V,E)$) consists of a finite set $V$ of vertices together with a set $E$ of unordered pairs of distinct vertices, called {\it edges}. Given a vertex $x$ of a graph $G=(V, E)$, the set $\{y\in V, \{x,y\}\in E\}$ is denoted by $N_{G}(x)$. With each subset $X$ of $V$ is associated the {\it subgraph} $G[X] = (X,E \cap \binom{X}{2})$ of $G$ induced by $X$. An isomorphism from a graph $G=(V,E)$ onto a graph $G'=(V',E')$ is a bijection $f$ from $V$ onto $V'$ such that for all $x$, $y \in V$, $\{x, y\} \in E$ if and only if $\{f(x), f(y)\} \in E'$.
We now introduce the graph $G_{2n}$ defined on $\mathbb{N}_{2n}$, where $n \geq 1$, as follows. For all $x$, $y \in \mathbb{N}_{2n}$, $\{x,y\} \in E(G_{2n})$ if and only if $\mid\! y-x \!\mid \geq n$ (see Figure 2).
\begin{center}
\psset{unit=0.8cm}
\begin{pspicture}(8, 5.4)
\psdots (1.9,4) (1.9,4.8) (1.9,0.8) (1.9,2.4) \put(1.85, 1.4) {$\vdots$} \put(1.85,3) {$\vdots$}
\put (0.9, 4.7){ \scriptsize$0$} \put (1, 3.9){\scriptsize$1$}  \put (1, 2.3){\scriptsize$i$}\put (0.7, 0.7){\scriptsize$n-1$}

\psdots (6.8,4) (6.8,4.8) (6.8,0.8)  (6.8,2.4) \put(6.75, 1.4) {$\vdots$} \put(6.75,3) {$\vdots$}
\put (7.5, 4.7){\scriptsize $n$} \put (7.3, 3.9){\scriptsize$n+1$} \put(7.3, 2.3){\scriptsize$i+n$}  \put (7.2, 0.7){\scriptsize$2n-1$}

\psline(1.9,4.8)(6.8,4.8) \psline(1.9,4.8)(6.8,4) \psline(1.9,4.8)(6.8,2.4) \psline(1.9,4.8)(6.8,0.8)

\psline(1.9,4)(6.8,4) \psline(1.9,4)(6.8,2.4) \psline(1.9,4)(6.8,0.8)

\psline (1.9,2.4)(6.8,2.4) \psline (1.9,2.4)(6.8,0.8)
 \psline(1.9,0.8)(6.8,0.8)

 \put(1.3,2.8){\oval(2,5.3)}
\put(7.4,2.8){\oval(2,5.3)}

\put (3.2, -0.1) {{\bf Fig.2.} $G_{2n}$}

\end{pspicture}
\end{center}

\vspace{0.15cm}

A graph $G$ is {\it connected} if for all $x \neq y \in V(G)$, there is a sequence $x_{0}=x, \dots, x_{m}=y$ of vertices of $G$ such that for all $i \in \mathbb{N}_{m}$, $\{x_{i}, x_{i+1}\} \in E(G)$. For example, the graph $G_{2n}$ is connected. A {\it connected component} of a graph $G$ is a maximal subset $X$ of $V(G)$ (with respect to inclusion) such that $G[X]$ is connected. The set of the connected components of $G$ is a partition of $V(G)$, denoted by $\mathcal{C}(G)$. Let $T = (V,A)$ be an indecomposable tournament. With each subset $X$ of $V$ such that $\mid\! X \!\mid \geq 3$ and $T[X]$ is indecomposable, is associated its {\it outside} graph $G^{T}_{X}$ defined by $V(G^{T}_{X})=V \setminus X$ and $E(G^{T}_{X}) = \{\{x,y\} \in \binom{V \setminus X}{2}: T[X \cup \{x,y\}]$ is indecomposable$\}$. We now present the characterization of partially critical tournaments.

\begin{thm}[\cite{MYS}] \label{partiellement critique}
Consider a tournament $T = (V,A)$ with a subset $X$ of $V$ such that $\mid\! X \!\mid \geq 3$ and $T[X]$ is indecomposable. The tournament $T$ is $T[X]$-critical if and only if the assertions below hold.
\begin{enumerate}
\item Ext$(X) = \emptyset$.
\item For all $u\in X$, the tournaments $T[X(u)\cup \{u\}]$ and $T[\langle X\rangle\cup \{u\}]$ are transitive.
\item For each $Q \in \mathcal{C}(G^{T}_{X})$, there is an isomorphism $f$ from $G_{2n}$ onto $G^{T}_{X}[Q]$ such that $Q_{1}$, $Q_{2}\in q^{T}_{X}$, where $Q_{1}=f(\mathbb{N}_{n})$ and $Q_{2}=f(\mathbb{N}_{2n}\setminus \mathbb{N}_{n})$. Moreover, for all $x\in Q_{i}\ (i=1\ or \ 2)$,
$\mid\! N_{G^{T}_{X}}(x) \!\mid= \mid\! N^{+}_{T[Q_{i}]}(x) \!\mid + 1$ $($resp. $n - \mid\! N^{+}_{T[Q_{i}]}(x) \!\mid)$ if $Q_{i} =X^{+}$ or $X^{-}(u)$ $($resp. $Q_{i}=X^{-}$ or $X^{+}(u))$, where~$u\in X$.
\end{enumerate}
\end{thm}

The next corollary follows from Theorem \ref{partiellement critique} and Lemma \ref{ER}.
\begin{cor}\label{000}
Let $T$ be a $T[X]$-critical tournament, $T$ is entirely determined up to isomorphy by giving $T[X]$, $q^{T}_{X}$ and $\mathcal{C}(G^{T}_{X})$. Moreover, the tournament $T$ is exactly determined by giving, in addition, either the graphs $G^{T}_{X}[Q]$ for any $Q \in \mathcal{C}(G^{T}_{X})$, or the transitive tournaments $T[Y]$ for any $Y\in q^{T}_{X}$.
\end{cor}

We underline the importance of Theorem \ref{partiellement critique} and Corollary \ref{000} in our description of the tournaments of the class $\mathcal{T}$. Indeed, these tournaments are introduced up to isomorphy as $C_{3}$-critical tournaments $T$ defined by giving $\mathcal{C}(G^{T}_{\mathbb{N}_{3}})$ in terms of the nonempty elements of $q^{T}_{\mathbb{N}_{3}}$. Figure 3 illustrates a tournament obtained from such information. We refer to [10, Discussion] for more details about this purpose.

We now introduce the class $\mathcal{H}$ (resp. $\mathcal{I}$, $\mathcal{J}$, $\mathcal{K}$, $\mathcal{L}$) of the $C_{3}$-critical tournaments $H$ (resp. $I$, $J$, $K$, $L$) such that:
\begin{itemize}
\item $\mathcal{C}(G^{H}_{\mathbb{N}_{3}})=\{\mathbb{N}_{3}^{+}(0)\cup \mathbb{N}_{3}^{-}, \mathbb{N}_{3}^{+}\cup \mathbb{N}_{3}^{-}(1)\}$ (see Figure~3);
\item $\mathcal{C}(G^{I}_{\mathbb{N}_{3}})=\{\mathbb{N}_{3}^{+}(0)\cup \mathbb{N}_{3}^{+}(2), \mathbb{N}_{3}^{+}(1)\cup \mathbb{N}_{3}^{-}(0)\}$;
\item $\mathcal{C}(G^{J}_{\mathbb{N}_{3}})=\{\mathbb{N}_{3}^{+}(1)\cup \mathbb{N}_{3}^{-}, \mathbb{N}_{3}^{-}(1)\cup \mathbb{N}_{3}^{-}(0)\}$;
\item $\mathcal{C}(G^{K}_{\mathbb{N}_{3}})=\{\mathbb{N}_{3}^{+}(1)\cup \mathbb{N}_{3}^{-}, \mathbb{N}_{3}^{+}(0)\cup \mathbb{N}_{3}^{-}(2)\}$;
\item $\mathcal{C}(G^{L}_{\mathbb{N}_{3}})=\{\mathbb{N}_{3}^{+}(1)\cup \mathbb{N}_{3}^{-}, \mathbb{N}_{3}^{+}(0)\cup \mathbb{N}_{3}^{-}(2), \mathbb{N}_{3}^{+} \cup \mathbb{N}_{3}^{-}(0)\}$.
\end{itemize}

Notice that for $\mathcal{X}=\mathcal{H}$, $\mathcal{I}$, $\mathcal{J}$ or $\mathcal{K}$, $\{\mid\!V(T)\!\mid:T\in \mathcal{X}\}=\{2n+1: n\geq 3\}$ and $\{\mid\!V(T)\!\mid:T\in \mathcal{L}\}=\{2n+1: n\geq 4\}$.
We denote by $\mathcal{H}^{\star}$ (resp. $\mathcal{I}^{\star}$, $\mathcal{J}^{\star}$, $\mathcal{K}^{\star}$, $\mathcal{L}^{\star}$) the class of the tournaments $T^{\star}$, where $T\in \mathcal{H}$ (resp. $\mathcal{I}$, $\mathcal{J}$, $\mathcal{K}$, $\mathcal{L}$).
\begin{rem}\label{dual}
We have $\mathcal{H}^{\star}=\mathcal{H}$ and $\mathcal{I}^{\star}=\mathcal{I}$.
\end{rem}
\bpr
Let $T\in \mathcal{H}$. The permutation $f$ of $V(T)$ defined by $f(1)=0$, $f(0)=1$ and $f(v)=v$ for all $v\in V(T)\setminus\{0,1\}$, is an isomorphism from $T^{\star}$ onto a tournament $T'$ of the class $\mathcal{H}$. Let now $T\in \mathcal{I}$ and let $x$ be the unique vertex of $\mathbb{N}_{3}^{+}(2)$ such that $\mid\!N^{+}_{T[\mathbb{N}_{3}^{+}(2)]}(x)\!\mid=0$. The permutation $g$ of $V(T)$ defined by $g(1)=0$, $g(0)=1$, $g(x)=2$, $g(2)=x$ and  $g(v)=v$ for $v\in V(T)\setminus\{0,1,2,x\}$, is an isomorphism from $T^{\star}$ onto a tournament $T'$ of the class $\mathcal{I}$.
\epr

By setting $\mathcal{M}=\mathcal{H}\cup \mathcal{I}\cup \mathcal{J} \cup \mathcal{J}^{\star} \cup \mathcal{K} \cup \mathcal{K}^{\star}\cup \mathcal{L}\cup \mathcal{L}^{\star}$, we state our main result as follows.
\begin{thm}\label{HIK2}
Up to isomorphy, the tournaments of the class $\mathcal{T}$ are those of the class $\mathcal{M}$. Moreover, for all $T\in\mathcal{M}$, we have  $V(T)\setminus W_{5}(T)=\sigma(T)=\{0,1\}$.
\end{thm}
\vspace{0.3 cm}

\begin{center}
\psset{unit=0.7cm}
\begin{pspicture}

\psline{->}(4,3.6)(5.9,3.6) \psline (5.8,3.6)(7.5,3.6)
\psline{->}(7.5,3.6)(6.53,2) \psline (6.58,2.08)(5.75,0.6)
\psline{->}(5.75,0.6)(4.73,2.29) \psline (4.77,2.23)(4,3.6)
\psdots (4,3.6) (7.5,3.6) (5.75,0.6)
\put (5.6, 0){2} \put (7.2, 3.8){1} \put (4.1, 3.8){$\tiny{0}$}

\psline{->}(4,3.6)(0.9,4.75)\psline(1,4.7)(-1.55,5.6)
\psframe[framearc=0.5](-3.25,5.5)(-1.25,9.5)
\put(-2.7,9.8 ){\small$\mathbb{N}_{3}^{+}(0)$}
\psdots (-1.75,9) (-1.75,8.5) (-1.75,6.5) (-1.75,6) \put(-1.82,7.3) {$\vdots$}
\put (-2.5, 8.9){$\scriptstyle{3}$} \put (-2.5, 8.4){$\scriptstyle{4}$} \put (-2.5, 6.4){$\scriptstyle{k}$} \put (-2.7, 5.9){$\scriptstyle{k+1}$}
\psline{->}(-1.75,8.5)(-1.75,8.88)\psline(-1.75,8.85) (-1.75,9) \psline{->}(-1.75,6)(-1.75,6.38)\psline(-1.75,6.35)(-1.75,6.5)
\psframe[framearc=0.5](0.75,5.5)(2.75,9.5)
\put(1.6,9.8 ){\small$\mathbb{N}_{3}^{-}$}
\psdots (1.25,9) (1.25,8.5) (1.25,6.5) (1.25,6) \put(1.18,7.3) {$\vdots$}
\put (1.75, 8.9){$\scriptstyle{k+2}$} \put (1.75, 8.4){$\scriptstyle{k+3}$} \put (1.6, 6.4){$\scriptstyle{2k-1}$} \put (1.85, 5.9){$\scriptstyle{2k}$}
\psline{->}(1.23,9)(1.23,8.62)\psline(1.23,8.65) (1.23,8.5) \psline{->}(1.23,6.5)(1.23,6.12)\psline(1.23,6.15)(1.23,6)
\psline{->}(3.05,7.65)(5.95,7.65)\psline(5.85,7.65)(8.45,7.65)
\psframe[framearc=0.5](-3.6,5.2)(3.05,10.3)
\psframe[framearc=0.5](8.45,5.2)(15.25,10.3)
\psline{->}(-1.75,9)(-0.15,9)\psline(-0.35,9)(1.25,9)
\psline{->}(-1.75,9)(-0.15,8.7)\psline(-0.2,8.7)(1.25,8.5)
\psline{->}(-1.75,9)(-0.15,7.7)\psline(-0.17,7.71)(1.25,6.5)
\psline{->}(-1.75,9)(-0.15,7.34)\psline(-0.2,7.37)(1.25,6)
\psline{->}(-1.75,8.5)(-0.15,8.5)\psline(-0.35,8.5)(1.25,8.5)
\psline{->}(-1.75,8.5)(-0.15,7.53)\psline(-0.2,7.54)(1.25,6.5)
\psline{->}(-1.75,8.5)(-0.15,7.18)\psline(-0.2,7.21)(1.25,6)
\psline{->}(-1.75,6.5) (-0.15,6.5)\psline(-0.35,6.5)(1.25,6.5)
\psline{->}(-1.75,6)(-0.15,6)\psline(-0.35,6)(1.25,6)
\psline{->}(-1.75,6.5)(-0.15,6.2)\psline(-0.2,6.2)(1.25,6)
\psline{->}(13.2,5.55)(10,4.5) \psline (10.08,4.53)(7.5,3.6)
\psframe[framearc=0.5](12.75,5.5)(14.9,9.5)
\put(13.3,9.8){\small$\mathbb{N}_{3}^{-}(1)$}
\psdots (13.25,9) (13.25,8.5) (13.25,6.5) (13.25,6) \put(13.18,7.3) {$\vdots$}
\put (13.5, 8.9){$\scriptstyle{n+k+1}$} \put (13.5, 8.4){$\scriptstyle{n+k+2}$} \put (13.65, 6.4){$\scriptstyle{2n-1}$} \put (13.85, 5.9){$\scriptstyle{2n}$}
\psline{->}(13.23,8.5)(13.23,8.88)\psline(13.23,8.85)(13.23,9) \psline{->}(13.23,6)(13.23,6.38) \psline(13.23,6.35)(13.23,6.5)
\psframe[framearc=0.5](8.75,5.5)(10.75,9.5)
\put(9.5,9.8){\small$\mathbb{N}_{3}^{+}$}
\psdots (10.25,9) (10.25,8.5) (10.25,6.5) (10.25,6) \put(10.18,7.3) {$\vdots$}
\put (9.1, 8.9){$\scriptstyle{2k+1}$} \put (9.1, 8.4){$\scriptstyle{2k+2}$} \put (8.9, 6.4){$\scriptstyle{n+k-1}$} \put (9.2, 5.9){$\scriptstyle{n+k}$}
\psline{->}(10.25,9) (10.25,8.62)\psline(10.25,8.65) (10.25,8.5) \psline{->}(10.25,6.5)(10.25,6.12)\psline(10.25,6.15)(10.25,6)
\psline{->}(10.25,9)(11.85,9)\psline(11.65,9) (13.25,9)
\psline{->}(10.25,9)(11.85,8.7)\psline(11.8,8.7)(13.25,8.5)
\psline{->}(10.25,9)(11.85,7.7) \psline(11.83,7.71)(13.25,6.5)
\psline{->}(10.25,9)(11.85,7.34)\psline(11.8,7.37)(13.25,6)
\psline{->}(10.25,8.5)(11.85,8.5)\psline(11.65,8.5)(13.25,8.5)
\psline{->}(10.25,8.5)(11.85,7.53)\psline(11.8,7.54)(13.25,6.5)
\psline{->}(10.25,8.5)(11.85,7.18)\psline(11.8,7.21)(13.25,6)
\psline{->}(10.25,6.5)(11.85,6.5)\psline(11.65,6.5)(13.25,6.5)
\psline{->}(10.25,6.5)(11.85,6.2)\psline(11.65,6.2)(13.25,6)
\psline{->}(10.25,6)(11.85,6)\psline(11.8,6)(13.25,6)

\psdots(-3.5,2.2)(-3.5,3)(-3.5,3.8)
\put(-3.3,3.7){{\scriptsize$T[\mathbb{N}_{3}^{-}]={\overrightarrow{\mathbb{N}_{2k+1}\setminus \mathbb{N}_{k+2}}}$;}}
\put(-3.3,2.9){{\scriptsize$T[\mathbb{N}_{3}^{+}(0)]={(\overrightarrow{\mathbb{N}_{k+2}\setminus \mathbb{N}_{3}})}^{\star}$;}}
\put(-3.3,2.1){{\scriptsize for all $(i,j)\in \mathbb{N}_{3}^{+}(0)\times\mathbb{N}_{3}^{-}$,}}
\put(-3.5,1.5){{\scriptsize $i\longrightarrow j$ if and only if $j-i\geq k-1$.}}
\psdots(9.5,2.2)(9.5,3)(9.5,3.8)
\put(9.7,3.7){{\scriptsize$T[\mathbb{N}_{3}^{+}]={\overrightarrow{\mathbb{N}_{n+k+1}\setminus \mathbb{N}_{2k+1}}}$;}}
\put(9.7,2.9){{\scriptsize$T[\mathbb{N}_{3}^{-}(1)]={(\overrightarrow{\mathbb{N}_{2n+1}\setminus \mathbb{N}_{n+k+1}})}^{\star}$;}}
\put(9.7,2.1){{\scriptsize for all $(i,j)\in \mathbb{N}_{3}^{+}\times\mathbb{N}_{3}^{-}(1)$,}}
\put(9.5,1.5){{\scriptsize $i\longrightarrow j$ if and only if $j-i\geq n-k$.}}

\put (0.3, -1) {{\bf Fig.3.}  A tournament $T$ of the class $\mathcal{H}$}
\end{pspicture}

\end{center}
\vspace{1 cm}

\section{Critical tournaments and tournaments omitting $W_{5}$}
\label{sec:2}
We begin by recalling the characterization of the critical tournaments and some of their properties. An indecomposable tournament $T=(V,A)$, with $\mid\! T\!\mid\geq3$, is {\it critical} if $\sigma(T)=\emptyset$, i.e., if all its vertices are critical.
In order to present the critical tournaments, characterized by J.H. Schmerl and W.T. Trotter in \cite{ST}, we introduce the tournaments $T_{2n+1}$ and $U_{2n+1}$ defined on $\mathbb{N}_{2n+1}$, where $n \geq 2$, as follows.
\begin{itemize}

\item $A(T_{2n+1}) = \{(i,j) : j-i \in \{1,\dots,n\}$ mod. $ 2n+1\}$ (see Figure 4).

\item $A(T_{2n+1})\setminus A(U_{2n+1})= A(T_{2n+1}[\{n+1,\dots,2n\}])$ (see Figure 5).

\end{itemize}
\vspace{0.5 cm}
\begin{center}
\psset{unit=0.85cm}
\begin{pspicture}(11,2)

\psdots (0,0) (1.5,0) (4.5,0) (6,0) (9,0) (10.5,0) (0.75,1.5) (5.25,1.5) (9.75,1.5)

\put (0.0,1.7){\scriptsize $n+1$} \put (4,1.7){\scriptsize $n+i+1$} \put (9.5,1.7){\scriptsize $2n$}

\put (-0.1, -0.4){\scriptsize ${0}$} \put (1.4, -0.4){\scriptsize $1$}

\put (4.26, -0.35) {$i$} \put (5.2, -0.35){\scriptsize $i+1$}

\put (8.55,-0.4){\scriptsize$n-1$} \put (10.35, -0.4){\scriptsize $n$}

\put(2.8,0) {$\ldots$}  \put (7.3,0){$\ldots$} \put (2.8,1.5){$\ldots$} \put (7.3,1.5){$\ldots$}

\psline{->}(1.5,0) (2.3,0)

\psline{->}(3.7,0) (4.5,0) \psline{->}(6,0) (6.8,0)

\psline{->}(8.2,0) (9,0)

\psline{->}(0,0)(0.85,0) \psline (0.75,0) (1.5,0)

\psline{->}(4.5,0)(5.35,0) \psline (5.3,0) (6,0)

\psline{->}(9,0)(9.85,0) \psline(9.8,0)(10.5,0)
\psline(0,0)(0.5,-0.4)\psline{->}(0.5,-0.4) (1,-0.4)
\psline(1.5,0)(2,-0.4)\psline{->}(2,-0.4) (2.5,-0.4)
\psline(4.5,0)(5,-0.4)\psline{->}(5,-0.4) (5.5,-0.4)
\psline(6,0)(6.5,-0.4)\psline{->}(6.5,-0.4) (7,-0.4)

\psline{->}(0.75,1.5) (1.6,1.5)
\psline{->}(4.4,1.5) (5.25,1.5)
\psline{->}(5.25,1.5) (6.1,1.5)
\psline{->}(8.9,1.5) (9.75,1.5)
\psline{->}(0.75,1.5) (0.348,0.63) \psline(0.348,0.63) (0,0)
\psline{->}(1.5,0)(1.1,0.75) \psline(1.1,0.75) (0.75,1.5)
\psline{->}(5.25,1.5) (4.848,0.63) \psline(4.848,0.63) (4.5,0)
\psline{->}(6,0)(5.598,0.75) \psline(5.598,0.75) (5.25,1.5)

\psline{->}(9.75,1.5) (9.348,0.63) \psline(9.348,0.63) (9,0)
\psline{->}(10.5,0)(10.098,0.75) \psline(10.098,0.75) (9.75,1.5)
\psline{->}(4.5,0)(4,0.3) \psline{->}(6,0)(5.5,0.3) \psline{->}(9,0)(8.5,0.3) \psline{->}(10.5,0)(10,0.3)
\psline{->}(9.75,1.5)(9.25,1.2) \psline{->}(5.25,1.5)(4.75,1.2)

\psline(0.75,1.5)(1.25,1.9)\psline{->}(1.25,1.9) (1.75,1.9)
\psline(5.25,1.5)(5.75,1.9)\psline{->}(5.75,1.9) (6.25,1.9)

\put (4, -1.2) {{\bf Fig.4.} $T_{2n+1}$}

\end{pspicture}

\end{center}

\vspace{2cm}

\begin{center}
\psset{unit=0.85cm}
\begin{pspicture}(11,1.7)

\psdots (0,0) (1.5,0) (4.5,0) (6,0) (9,0) (10.5,0) (0.75,1.5) (5.25,1.5) (9.75,1.5)

\put (0.0,1.7){\scriptsize $n+1$} \put (5.1,1.7){\scriptsize $n+i+1$} \put (9.5,1.7){\scriptsize $2n$}

\put (-0.1, -0.4){\scriptsize ${0}$} \put (1.4, -0.4){\scriptsize $1$}

\put (4.26, -0.35) {\scriptsize $i$} \put (5.2, -0.35){\scriptsize $i+1$}

\put (8.55,-0.4){\scriptsize $n-1$} \put (10.35, -0.4){\scriptsize $n$}

\put(2.8,0) {$\ldots$}  \put (7.3,0){$\ldots$} \put (2.8,1.5){$\ldots$} \put (7.3,1.5){$\ldots$}

\psline{->}(1.5,0) (2.3,0)

\psline{->}(3.7,0) (4.5,0) \psline{->}(6,0) (6.8,0)

\psline{->}(8.2,0) (9,0)

\psline{->}(0,0)(0.85,0) \psline (0.75,0) (1.5,0)

\psline{->}(4.5,0)(5.35,0) \psline (5.3,0) (6,0)

\psline{->}(9,0)(9.85,0) \psline(9.8,0)(10.5,0)

\psline(0,0)(0.5,-0.4)\psline{->}(0.5,-0.4) (1,-0.4)
\psline(1.5,0)(2,-0.4)\psline{->}(2,-0.4) (2.5,-0.4)
\psline(4.5,0)(5,-0.4)\psline{->}(5,-0.4) (5.5,-0.4)
\psline(6,0)(6.5,-0.4)\psline{->}(6.5,-0.4) (7,-0.4)

\psline{<-}(0.75,1.5) (1.6,1.5)
\psline{<-}(4.4,1.5) (5.25,1.5)
\psline{<-}(5.25,1.5) (6.1,1.5)
\psline{<-}(8.9,1.5) (9.75,1.5)
\psline{->}(0.75,1.5) (0.348,0.63) \psline(0.348,0.63) (0,0)
\psline{->}(1.5,0)(1.1,0.75) \psline(1.1,0.75) (0.75,1.5)
\psline{->}(5.25,1.5) (4.848,0.63) \psline(4.848,0.63) (4.5,0)
\psline{->}(6,0)(5.598,0.75) \psline(5.598,0.75) (5.25,1.5)

\psline{->}(9.75,1.5) (9.348,0.63) \psline(9.348,0.63) (9,0)
\psline{->}(10.5,0)(10.098,0.75) \psline(10.098,0.75) (9.75,1.5)
\psline{->}(4.5,0)(4,0.3) \psline{->}(6,0)(5.5,0.3) \psline{->}(9,0)(8.5,0.3) \psline{->}(10.5,0)(10,0.3)
\psline{->}(9.75,1.5)(9.25,1.2) \psline{->}(5.25,1.5)(4.75,1.2)

\psline(9.75,1.5)(9.25,1.9)\psline{->}(9.25,1.9) (8.75,1.9)
\psline(5.25,1.5)(4.75,1.9)\psline{->}(4.75,1.9) (4.25,1.9)

\put (4, -1.2) {{\bf Fig.5.} $U_{2n+1}$}
\end{pspicture}

\end{center}

\vspace{1.2cm}

\begin{thm} [\cite{ST}] \label{P1}
Up to isomorphy, $T_{2n+1}$, $U_{2n+1}$ and $W_{2n+1}$, where $n \geq 2$, are the only critical tournaments.
\end{thm}

Notice that a critical tournament is isomorphic to its dual. Moreover, as a tournament on 4 vertices is decomposable, we have:
\begin{Fact}\label{indec a 5sommets}
Up to isomorphy, $T_{5}$, $U_{5}$ and $W_{5}$ are the only indecomposable tournaments on 5 vertices.
\end{Fact}

The next remark follows from the definition of the critical tournaments.
\begin{rem}[\cite{YJP}] \label{abrit entre critique}
Up to isomorphy, the indecomposable subtournaments on at least 5 vertices of $T_{2n+1}$ $($resp. $U_{2n+1}$, $W_{2n+1})$, where $n \geq 2$,
are the tournaments $T_{2m+1}$ $($resp. $U_{2m+1}$, $W_{2m+1})$, where  $2 \leq m \leq n$.
\end{rem}

To recall the characterization of the indecomposable tournaments omitting $W_{5}$, we introduce the {\it Paley} tournament $P_{7}$ defined on $\mathbb{N}_{7}$
by $A(P_{7}) = \{(i,j) : j-i \in \{1, 2, 4\}$ mod.  $7\}$.
Notice that for all $x \neq y \in \mathbb{N}_{7}$, $P_{7}-x \simeq P_{7}-y$, and let $B_{6} = P_{7}-6$.
\begin{thm} [\cite{BJL}] \label{Latka}
Up to isomorphy, the indecomposable tournaments on at least $5$ vertices and omitting $W_{5}$
are the tournaments $B_{6}$, $P_{7}$, $T_{2n+1}$ and $U_{2n+1}$, where $n \geq 2$.
\end{thm}
\section{Some useful configurations}
\label{sec:3}
In this section, we introduce number of configurations which occur in the proof of Theorem \ref{HIK2}. These configurations involve mainly partially critical tournaments. We begin by the two following lemmas obtained in \cite{HIK'}.
\begin{lem}[\cite{HIK'}] \label{B6}
If $B_{6}$ embeds into an indecomposable tournament $T$ on 7 vertices and if $T\not\simeq P_{7}$, then $\mid\!W_{5}(T)\!\mid =7$.
\end{lem}
\begin{sloppypar}
\begin{lem} [\cite{HIK'}] \label{lemme U5}
Let $T$ be a $U_{5}$-critical tournament on 7 vertices. If $T \not \simeq U_{7}$, then $W_{5}(T)~\cap~\{3, 4\}\neq~\emptyset$.
\end{lem}
\end{sloppypar}
Lemma \ref{une composante} specifies the $C_{3}$-critical tournaments with a connected outside graph. It follows from the examination of the different possible configurations obtained by using Theorem \ref{partiellement critique}.

\begin{lem} \label{une composante}
Given a $C_{3}$-critical tournament $T$ on at least 5 vertices, if $G^{T}_{\mathbb{N}_{3}}$ is connected, then $T$ is critical. More precisely, the different configurations are as follows where $i\in \mathbb{N}_{3}$ and $i+1$ is considered modulo 3.
\begin{enumerate}
\item If $\mathcal{C}(G^{T}_{\mathbb{N}_{3}})=\{\mathbb{N}_{3}^{-}(i)\cup \mathbb{N}_{3}^{+}(i+1)\}$, then $T\simeq T_{2n+1}$ for some $n\geq2$.

\item If $\mathcal{C}(G^{T}_{\mathbb{N}_{3}})=\{\mathbb{N}_{3}^{-}\cup \mathbb{N}_{3}^{+}(i)\}$, $\{\mathbb{N}_{3}^{+}\cup \mathbb{N}_{3}^{-}(i)\}$, $\{\mathbb{N}_{3}^{+}(i)\cup \mathbb{N}_{3}^{+}(i+1)\}$ or $\{\mathbb{N}_{3}^{-}(i)\cup \mathbb{N}_{3}^{-}(i+1)\}$, then $T\simeq U_{2n+1}$ for some $n\geq2$.

\item If $\mathcal{C}(G^{T}_{\mathbb{N}_{3}})=\{\mathbb{N}_{3}^{-}\cup \mathbb{N}_{3}^{-}(i)\}$, $\{\mathbb{N}_{3}^{+}\cup \mathbb{N}_{3}^{+}(i)\}$ or $\{\mathbb{N}_{3}^{+}(i)\cup \mathbb{N}_{3}^{-}(i+1)\}$, then $T\simeq W_{2n+1}$ for some $n\geq2$.

\end{enumerate}
\end{lem}

For a transitive tournament $T$, recall that $\min T$ denotes its smallest element and
$\max T$ its largest.

\begin{lem}\label{T-2n+1}
Given a $C_{3}$-critical tournament $T$ on at least 5 vertices, if $T[\mathbb{N}_{3}\cup e]\simeq T_{5}$ for all $e \in E(G^{T}_{\mathbb{N}_{3}})$, then $T\simeq T_{2n+1}$ for some $n\geq 2$.
\end{lem}
\bpr
Let $T$ be a $C_{3}$-critical tournament on at least 5 vertices such that for all $e \in E(G^{T}_{\mathbb{N}_{3}})$, $T[\mathbb{N}_{3}\cup e]\simeq T_{5}$.
Given $e \in E(G^{T}_{\mathbb{N}_{3}})$, by using Lemma \ref{une composante} and Remark \ref{abrit entre critique}, $e=\{v, v'\}$ where
$v\in \mathbb{N}_{3}^{-}(i)$, $v'\in \mathbb{N}_{3}^{+}(i+1)$, $i\in \mathbb{N}_{3}$ and $i+1$ is considered modulo 3. Then, by Theorem \ref{partiellement critique}, the connected components of $T$ are the nonempty elements of the family $\{\mathbb{N}_{3}^{-}(j) \cup \mathbb{N}_{3}^{+}(j+1)\}_{j\in \mathbb{N}_{3}}$, where $j+1$ is considered modulo 3. The tournament $T$ is critical. Indeed, by using Theorem \ref{partiellement critique}, for each $k \in \mathbb{N}_{3}$, $\{\max T[\mathbb{N}_{3}^{+}(k+1)\cup \{k+1\}], \min T[\mathbb{N}_{3}^{-}(k+2)\cup \{k+2\}]\}$, where $k+1$ and $k+2$ are considered modulo 3, is a non-trivial interval of $T-k$. It follows that $T\simeq T_{2n+1}$ for some $n\geq 2$ by Remark~\ref{abrit entre critique}.
\epr

\begin{lem}\label{U-2n+1}
Given a $U_{5}$-critical tournament, if $T[\mathbb{N}_{5}\cup e]\simeq U_{7}$ for all $e \in E(G^{T}_{\mathbb{N}_{5}})$, then $T\simeq U_{2n+1}$ for some $n\geq 2$.
\end{lem}
\begin{sloppypar}
\bpr
The subsets $X$ of $\mathbb{N}_{7}$ such that $U_{7}[X] \simeq U_{5}$ are the sets $\mathbb{N}_{7} \setminus \{i,j\}$, where $\{i,j\} = \{0,4\}$, $\{4,1\}$, $\{1,5\}$, $\{5,2\}$, $\{2,6\}$ or $\{6,3\}$. By observing $q^{U_{7}}_{X}$ for such subsets $X$ and by Theorem \ref{partiellement critique}, we deduce that the elements of $\mathcal{C}(G^{T}_{\mathbb{N}_{5}})$ are the nonempty elements among the following six sets: $\mathbb{N}_{5}^{+} \cup \mathbb{N}_{5}^{-}(0)$,  $\mathbb{N}_{5}^{+}(0) \cup \mathbb{N}_{5}^{+}(3)$, $\mathbb{N}_{5} ^{-}(1)\cup \mathbb{N}_{5}^{-}(3)$, $\mathbb{N}_{5}^{+}(1) \cup \mathbb{N}_{5}^{+}(4)$,  $\mathbb{N}_{5}^{-}(2) \cup \mathbb{N}_{5}^{-}(4)$ and $\mathbb{N}_{5}^{-} \cup \mathbb{N}_{5}^{+}(2)$. Suppose first that $\mid\!\mathcal{C}(G^{T}_{\mathbb{N}_{5}})\!\mid =6$. The tournament $T$ is critical. Indeed, by using Theorem \ref{partiellement critique}, $\{\min T[ \mathbb{N}_{5}^{+}], \max T[ \mathbb{N}_{5}^{+}(3)]\}$ (resp. $\{ \min T[ \mathbb{N}_{5}^{-}(3)], \max T[ \mathbb{N}_{5}^{+}(4)]\}$,  $\{ \min T[ \mathbb{N}_{5}^{-}(4)], \max T[ \mathbb{N}_{5}^{-}] \}$, $\{ \min T[ \mathbb{N}_{5}^{-}(1)], \max T[ \mathbb{N}_{5}^{+}(0)] \}$,  $\{ \min T[ \mathbb{N}_{5}^{-}(2)], \max T[ \mathbb{N}_{5}^{+}(1)] \}$) is a non-trivial interval of $T-0$ (resp. $T-1$, $T-2$, $T-3$, $T-4$). By Remark \ref{abrit entre critique}, $T\simeq U_{2n+1}$ for some $n\geq 8$. Suppose now that $\mid\!\mathcal{C}(G^{T}_{\mathbb{N}_{5}})\!\mid \leq 5$. Then, $T$ embeds into a $U_{5}$-critical tournament $T'$ with $\mid\!\mathcal{C}(G^{T'}_{\mathbb{N}_{5}})\!\mid =6$. By the first case, $T'\simeq U_{2n+1}$ for some $n\geq 8$ and thus  $T\simeq U_{2n+1}$ for some $n\geq 2$ by Remark \ref{abrit entre critique}. \epr
\end{sloppypar}

\begin{lem}\label{G-2n-2}
 Let $T=(V,A)$ be a $T[X]$-critical tournament with $\mid\!V\setminus X\!\mid\geq$~$2$, let $Q=\mathbb{N}_{2n}$ be a connected component of $G^{T}_{X}$ such that $G^{T}_{X}[Q]=G_{2n}$ and let $e=\{i, i+n\}$, where $i \in \mathbb{N}_{n}$. Then, the tournament $T-e$ is  $T[X]$-critical. Moreover, $Q$ is included in any subset $Z$ of $V$ such that $T[Z]\simeq W_{5}$ and $Z\cap (V\setminus (Q\cup W_{5}(T-e))\neq \emptyset$.
\end{lem}
\bpr
For $n\geq 2$, the function
\begin{eqnarray*}
  & f_{i} : & Q \setminus e \longrightarrow \mathbb{N}_{2n-2} \\ \\
  \nonumber & & \ \ \ k \ \ \   \longmapsto \left\{
                       \begin{array}{ll}
                         \ \ \  k \, & \hbox{ if \ $0 \leq k \leq i-1$} \\
                         k-1 \, & \hbox{ if \ $i+1\leq k \leq n+i-1$} \\
                         k-2 \, & \hbox{ if \ $n+i+1\leq k \leq 2n-1$,}
                       \end{array}
                     \right. \\
\end{eqnarray*}
is an isomorphism from $G_{2n}-e$ onto $G_{2n-2}$. It follows from Theorem \ref{partiellement critique} that $T-e$ is $T[X]$-critical. Now, suppose that there is $Z \subseteq V$ such that $T[Z]\simeq W_{5}$ and $Z \cap (V\setminus (Q\cup W_{5}(T-e))\neq \emptyset$. So we have $\mid \! Z\cap e \!\mid=1 $ or $e \subset Z$. Suppose for a contradiction that $\mid \! Z\cap e \!\mid=1$, and set $\{z\}=Z\cap e$. As Ext$(V\setminus e)=\emptyset$, then by Lemma \ref{ER}, either $z\in \langle V'\rangle$ or $z\in V'(u)$, where $V'=V\setminus e$ and $u\in V'$. If $z\in \langle V'\rangle$, then $Z\setminus \{z\}$ is a non-trivial interval of $T[Z]$, a contradiction. If $z\in V'(u)$, then $u\not\in Z$, otherwise $\{u,z\}$ is a non-trivial interval of $T[Z]$. Thus, $T[Z']\simeq W_{5}$, where $Z'=(Z\setminus \{z\})\cup\{u\}\subset V\setminus e$. A contradiction because $Z'\cap (V\setminus W_{5}(T-e))\neq \emptyset$. Finally, for all $e' \in \{ \{j, j+n\}:  j \in \mathbb{N}_{n}\}$, the bijection $f$ from $V\setminus e$ onto $V \setminus e'$, defined by
$\restriction{f}{V\setminus Q}$ $= \text{Id}_{V \setminus Q}$ and $\restriction{f}{Q \setminus e}$ $= f_{j}^{-1} \circ f_{i}$, is an isomorphism from $T - e$ onto $T - e'$. It follows that $V\setminus( Q\cup W_{5}(T-e')) = V\setminus( Q\cup W_{5}(T-e))$. Thus, as proved above, $e' \subset Z$, so that $Q \subset Z$.
\epr

\section{Proof of Theorem \ref{HIK2}}
\label{sec:4}

We begin by establishing the partial criticality structure of the tournaments of the class $\mathcal{T}$. For this purpose, we use the notion of minimal tournaments for two vertices. Given an indecomposable tournament $T = (V,A)$ of cardinality $\geq 3$ and two distinct vertices $x \neq y \in V$,
$T$ is said to be {\it minimal} for $\{x,y\}$ (or {\it $\{x,y\}$-minimal}) when for all
proper subset $X$ of $V$, if $\{x,y\} \subset X$ ($\mid\! X \!\mid \geq 3$),
then $T[X]$ is decomposable. These tournaments were introduced and characterized  by A. Cournier and P. Ille in \cite{CI}.
 From this characterization, the following fact, observed in \cite{HIK}, is obtained by a simple and quick verification.
\begin{Fact} [\cite{HIK,CI}] \label{Minimaux}
Up to isomorphy, the tournaments $C_{3}$ and $U_{5}$ are the unique minimal tournaments for two vertices $T$ such that $\mid\! W_{5}(T)\!\mid \leq \mid\! T\!\mid-2$. Moreover, $\{3,4\}$ is the unique unordered pair of vertices for which $U_{5}$ is minimal.
\end{Fact}
\begin{prop}\label{3cycle}
Let $T=(V,A)$ be a tournament of the class $\mathcal{T}$. Then, the vertices of $W_{5}(T)$ are critical and there exists $z\in W_{5}(T)$ such that $T[(V\setminus W_{5}(T))\cup\{z\}]\simeq C_{3}$. In particular, $T$ is $T[(V\setminus W_{5}(T))\cup\{z\}]$-critical.
\end{prop}
\bpr
By Theorem \ref{HIK}, $\mid\! T \!\mid$ is odd and $\geq 7$. First, suppose by contradiction, that there is $\alpha \in W_{5}(T)$ such that $T- \alpha$ is indecomposable. Since $\mid\!T - \alpha\!\mid$ is even and $\geq 6$ with $\mid\! V(T-\alpha)\setminus W_{5}(T - \alpha) \!\mid \geq 2$, then by Theorems \ref{HIK} and \ref{Latka}, $T-\alpha \simeq B_{6}$ and $T \not\simeq P_{7}$. A contradiction by Lemma \ref{B6}. Second, let $X$ be a minimal subset of $V$ such that $V\setminus W_{5}(T)\subset X$ ($\mid\!X\!\mid\geq 3 $) and $T[X]$ is indecomposable, so that $T[X]$ is $(V\setminus W_{5}(T))$-minimal. By Fact \ref{Minimaux}, $T[X]\simeq C_{3}$ or $U_{5}$. Suppose, toward a contradiction that $T[X]\simeq U_{5}$ and take $T[X]=U_{5}$. By Fact \ref{Minimaux}, $V\setminus W_{5}(T)=\{3, 4\}$. As $T$ is $U_{5}$-critical, then by Lemma \ref{U-2n+1} and Theorem \ref{Latka}, there exists $e\in E(G^{T}_{X})$ such that $T[X\cup e]$ is indecomposable and not isomorphic to $U_{7}$. It follows from Lemma \ref{lemme U5}, that there exists a subset $Z$ of $X\cup e$ such that $T[Z]\simeq W_{5}$ and $Z\cap (V\setminus W_{5}(T))\neq \emptyset$, a contradiction.
\epr

Now, we prove Theorem \ref{HIK2} for tournaments on 7 vertices.
\begin{prop}\label{7 sommets}
Up to isomorphy, the class $\mathcal{M}$ and the class $\mathcal{T}$ have the same tournaments on 7 vertices. Moreover, for each tournament $T$ on 7 vertices of the class $\mathcal{M}$, we have  $V(T)\setminus W_{5}(T)=\sigma(T)=\{0,1\}$.
\end{prop}
 \bpr
Let $T=(V,A)$ be a tournament on 7 vertices of the class $\mathcal{M}$. $T\in \mathcal{M}\setminus (\mathcal{L}\cup \mathcal{L}^{\star})$ because the tournaments of the class $\mathcal{L}$ have at least 9 vertices. Let $e\in E(G^{T}_{\mathbb{N}_{3}})$. By Lemma~\ref{une composante}, $T-e\simeq U_{5}$ or $T_{5}$. By Lemma \ref{G-2n-2}, if there exists a subset $Z\subset V$ such that $T[Z]\simeq W_{5}$, then $e \subset Z$. It follows that $V\setminus \mathbb{N}_{3} \subset Z$. Thus $V \setminus W_{5}(T)=\{0,1\}$ by verifying that $T-\{1,2\}\not\simeq W_{5}$, $T-\{0,2\}\not \simeq W_{5}$ and $T-\{0,1\}\simeq W_{5}$. As $T$ is $C_{3}$-critical, $\sigma(T)=\{0,1\}$ from the following. First, $T-2$ is decomposable because $\{0\}\cup \mathbb{N}_{3}^{-}\cup \mathbb{N}_{3}^{+}(0)$ (resp. $\{1\}\cup \mathbb{N}_{3}^{+}(0)$, $\{0,1\}\cup \mathbb{N}_{3}^{-}(0)\cup \mathbb{N}_{3}^{-}(1)$, $\{1\}\cup \mathbb{N}_{3}^{+}(0)$) is a non-trivial interval of $T-2$ if $T\in \mathcal{H}$ (resp. $\mathcal{I}$, $\mathcal{J}$, $\mathcal{K}$). Second, by Lemma \ref{ER}, we have Ext$(X) = \{0, 1\}$, where $X= V\setminus \{0, 1\}$, because $\{0,1\}\cap\langle X\rangle=\emptyset$, and for all $u\in X$, $\{0,1\}\cap X(u)=\emptyset$ because $V\setminus W_{5}(T)=\{0,1\}$.

Conversely, let $T$ be a tournament on 7 vertices of the class $\mathcal{T}$. By Proposition \ref{3cycle}, we can assume that $T$ is $C_{3}$-critical with $V(T)\setminus W_{5}(T)\subset \mathbb{N}_{3}$. By Lemma \ref{une composante} and Theorem \ref{Latka}, $\mid\!\mathcal{C}(G^{T}_{\mathbb{N}_{3}})\!\mid=2$. We distinguish the following cases.
\begin{itemize}
\item[$\bullet$] $\mathbb{N}_{3}^{+}\neq \emptyset$ and $\mathbb{N}_{3}^{-}\neq \emptyset$. By Theorem \ref{partiellement critique}, $\mid\! \mathbb{N}_{3}^{-} \!\mid = \mid\! \mathbb{N}_{3}^{+} \!\mid =1$. So, we can assume that $\mathbb{N}_{3}(0)\neq \emptyset$ and $\mathbb{N}_{3}(2)=\emptyset$. It suffices to verify that $\mid\! \mathbb{N}_{3}(0) \!\mid = \mid\! \mathbb{N}_{3}^{+}(0) \!\mid =1$ because, in this case, by using Theorem \ref{partiellement critique} and Lemma \ref{une composante}, $T\in \mathcal{H}$. By using again Theorem~\ref{partiellement critique} and Lemma \ref{une composante}, we verify the following. First, if $\mid\! \mathbb{N}_{3}(0) \!\mid =2$, then $\mathcal{C}(G^{T}_{\mathbb{N}_{3}})=\{\mathbb{N}_{3}^{+}\cup \mathbb{N}_{3}^{-}(0), \mathbb{N}_{3}^{-}\cup \mathbb{N}_{3}^{+} (0) \}$. Therefore, $T- \{0,1\}\simeq T- \{0,2\}\simeq W_{5}$, a contradiction. Second, if $\mid\! \mathbb{N}_{3}^{-}(0) \!\mid =1$, then $\mathcal{C}(G^{T}_{\mathbb{N}_{3}})=\{\mathbb{N}_{3}^{+}\cup \mathbb{N}_{3}^{-}(0), \mathbb{N}_{3}^{-}\cup \mathbb{N}_{3}^{+} (1) \}$. Therefore, $T\simeq U_{7}$, a contradiction by Theorem \ref{Latka}.
\begin{sloppypar}
\item[$\bullet$] $\langle\mathbb{N}_{3}\rangle= \emptyset$. By Theorem \ref{partiellement critique}, we can assume that $\mid\! \mathbb{N}_{3}^{-}(0) \!\mid = \mid\! \mathbb{N}_{3}^{+}(0) \!\mid=1$. We have $\mid\!~\!\mathbb{N}_{3}(1)\!\mid=1$. Otherwise, by Theorem \ref{partiellement critique} and Lemma \ref{une composante}, we can suppose that $\mathcal{C}(G^{T}_{\mathbb{N}_{3}})=\{\mathbb{N}_{3}^{+}(0)\cup \mathbb{N}_{3}^{+}(1), \mathbb{N}_{3}^{-}(0)\cup \mathbb{N}_{3}^{-}(1) \}$. Therefore, $T-\{1,2\}\simeq T-\{0,2\}\simeq W_{5}$, a contradiction.
We have also $\mathcal{C}(G^{T}_{\mathbb{N}_{3}})=\{\mathbb{N}_{3}^{+}(0)\cup \mathbb{N}_{3}(2), \mathbb{N}_{3}^{-}(0)\cup \mathbb{N}_{3}(1) \}$. Otherwise, again by Theorem \ref{partiellement critique} and Lemma \ref{une composante}, $\mathcal{C}(G^{T}_{\mathbb{N}_{3}})=\{\mathbb{N}_{3}^{+}(0)\cup \mathbb{N}_{3}^{+}(1), \mathbb{N}_{3}^{-}(0)\cup \mathbb{N}_{3}^{-}(2) \}$, so that $T\simeq U_{7}$, a contradiction by Theorem \ref{Latka}.
So, we distinguish four cases. If $\mid\!\mathbb{N}_{3}^{-}(2)\!\mid=\mid\!\mathbb{N}_{3}^{+}(1)\!\mid=1$, then $T\simeq T_{7}$, which contradicts Theorem \ref{Latka}.
If  $\mid\!\mathbb{N}_{3}^{+}(2)\!\mid=\mid\!\mathbb{N}_{3}^{-}(1)\!\mid=1$, then $T-\{0,2\}\simeq T-\{0,1\} \simeq W_{5}$, a contradiction. If $\mid\!\mathbb{N}_{3}^{+}(2)\!\mid=\mid\!\mathbb{N}_{3}^{+}(1)\!\mid=1$, then $T\in \mathcal{I}$. If $\mid\!\mathbb{N}_{3}^{-}(2)\!\mid=\mid\!\mathbb{N}_{3}^{-}(1)\!\mid=1$, then $T$ is isomorphic to a tournament of the class $\mathcal{I}$ with $V(T)\setminus W_{5}(T)=\{0,2\}$.
\end{sloppypar}
\item[$\bullet$] $\emptyset\neq \langle\mathbb{N}_{3}\rangle\in q^{T}_{\mathbb{N}_{3}}$. By interchanging $T$ and $T^{\star}$, we can suppose that $\langle\mathbb{N}_{3}\rangle = \mathbb{N}_{3}^{-}$. In this case, $\mid\!\mathbb{N}_{3}^{-}\!\mid=1$ by Theorem \ref{partiellement critique}.
First, suppose that $\mid\!\mathbb{N}_{3}(0)\!\mid=2$ and $\mid\!\mathbb{N}_{3}(1)\!\mid=1$. By Theorem \ref{partiellement critique} and Lemma \ref{une composante}, $\mathcal{C}(G^{T}_{\mathbb{N}_{3}})=\{\mathbb{N}_{3}^{+}(0)\cup \mathbb{N}_{3}^{-}, \mathbb{N}_{3}^{-}(0)\cup \mathbb{N}_{3}(1)\}$. We have $\mid\!\mathbb{N}_{3}^{+}(1)\!\mid=1$, otherwise $T\simeq U_{7}$, a contradiction by Theorem \ref{Latka}. Thus, $T$ is isomorphic to a tournament of the class $\mathcal{K}$ with $V(T)\setminus W_{5}(T)=\{0,2\}$. Second, suppose that $\mid\!\mathbb{N}_{3}(0)\!\mid=1$ and $\mid\!\mathbb{N}_{3}(1)\!\mid=2$. Again by Theorem \ref{partiellement critique} and Lemma~\ref{une composante}, $\mathcal{C}(G^{T}_{\mathbb{N}_{3}})=\{\mathbb{N}_{3}^{+}(1)\cup \mathbb{N}_{3}^{-}, \mathbb{N}_{3}^{-}(1)\cup \mathbb{N}_{3}^{-}(0)\}$, so that $T\in \mathcal{J}$. Lastly, suppose that $\mid~\!\mathbb{N}_{3}(0)\!\mid=\mid\!~\mathbb{N}_{3}(1)\!\mid=1$. By Theorem \ref{partiellement critique} and Lemma \ref{une composante}, we can suppose that $\mathcal{C}(G^{T}_{\mathbb{N}_{3}})=\{\mathbb{N}_{3}^{+}(1)\cup \mathbb{N}_{3}^{-}, \mathbb{N}_{3}(0)\cup \mathbb{N}_{3}(2)\}$. By Lemma \ref{une composante}, we distinguish only three cases. If $\mid\!\mathbb{N}_{3}^{-}(2)\!\mid=\mid\!\mathbb{N}_{3}^{-}(0)\!\mid=1$, then $T- \{0,1\}\simeq T- \{1,2\}\simeq W_{5}$, a contradiction. If $\mid\!\mathbb{N}_{3}^{+}(0)\!\mid=\mid\!\mathbb{N}_{3}^{+}(2)\!\mid=1$, then $T\simeq U_{7}$, which contradicts Theorem \ref{Latka}. If $\mid\!~\mathbb{N}_{3}^{-}(2)\!\mid=\mid\!\mathbb{N}_{3}^{+}(0)\!\mid=1$, then $T\in \mathcal{K}$. \epr

\end{itemize}

We complete our structural study of the tournaments of the class $\mathcal{T}$ by the following two corollaries.

\begin{cor} \label{abrit 7 sommets}

Let $T$ be a $C_{3}$-critical tournament such that $V(T)\setminus W_{5}(T)=\{0,1\}$. Then, there exist $Q\neq Q'\in \mathcal{C}(G^{T}_{\mathbb{N}_{3}})$ and a tournament $R$ on 7 vertices of the class $\mathcal{M}$ such that for all $e \in E(G^{T}_{\mathbb{N}_{3}}[Q])$ and for all $e'\in E(G^{T}_{\mathbb{N}_{3}}[Q'])$, there exists an isomorphism $f$ from $R$ onto $T[\mathbb{N}_{3}\cup e\cup e']$. Moreover, $f(0)=0$, $f(1)=1$ and we have:
\begin{enumerate}
\item If $R\in \mathcal{H}\cup \mathcal{J} \cup \mathcal{J}^{\star}$, then $f(2)=2$;
\item If $R\in \mathcal{I}\cup \mathcal{K} \cup \mathcal{K}^{\star}$, then $f(2)=2$ or $ \mathbb{N}_{3}(2)=\{f(2)\}$.
\end{enumerate}
\end{cor}
\bpr
To begin, notice the following remark: given a $D[X]$-critical tournament $D$, for any edges $a$ and $b$ belonging to a same connected component of $G^{D}_{X}$, we have $D[X \cup a]\simeq D[X \cup b]$. So, by Fact \ref{indec a 5sommets}, Lemma \ref{T-2n+1} and Theorem \ref{Latka}, there exists $Q \in \mathcal{C}(G^{T}_{\mathbb{N}_{3}})$ such that for all $a\in E(G^{T}_{\mathbb{N}_{3}}[Q])$, $T[\mathbb{N}_{3} \cup a] \simeq U_{5}$. By Lemma \ref{une composante} and Remark~\ref{abrit entre critique}, the tournament $T[\mathbb{N}_{3}\cup Q]$ is isomorphic to $U_{2n+1}$, for some $n\geq 2$, and does not admit an indecomposable subtournament on 7 vertices other than $U_{7}$. So, by Lemma \ref{U-2n+1}, Theorem \ref{Latka} and the remark above, there exists $Q'\in \mathcal{C}(G^{T}_{\mathbb{N}_{3}})\setminus \{Q\} $ such that for all $e \in E(G^{T}_{\mathbb{N}_{3}}[Q])$ and for all $e' \in E(G^{T}_{\mathbb{N}_{3}}[Q'])$, $T[\mathbb{N}_{3}\cup e\cup e']$ is indecomposable and not isomorphic to $U_{7}$. Moreover, $T[\mathbb{N}_{3} \cup e\cup e'] \not \simeq P_{7}$ because the vertices of $P_{7}$ are all non-critical. Likewise, $T[\mathbb{N}_{3} \cup e \cup e'] \not \simeq T_{7}$ by Remark \ref{abrit entre critique}. It follows from Theorem \ref{Latka} and Proposition \ref{7 sommets}, that there exists an isomorphism $f$ from a tournament $R$ on 7 vertices of the class $\mathcal{M}$ onto $T[\mathbb{N}_{3} \cup e \cup e']$. As $(0,1)\in A(R)\cap A(T)$ and $V(R)\setminus W_{5}(R)=V(T)\setminus W_{5}(T)=\{0,1\}$ by Proposition \ref{7 sommets}, then $f$ fixes 0 and 1. If $R\in \mathcal{H}\cup \mathcal{J} \cup \mathcal{J}^{\star}$, then $f$ fixes 2 because 2 is the unique vertex $x$ of $R$ such that $R[\{0,1,x\}]\simeq C_{3}$. If $R\in \mathcal{I}\cup \mathcal{K} \cup \mathcal{K}^{\star}$, then $\mid\!\{x\in V(R) :  R[\{0,1,x\}]\simeq C_{3}\}\!\mid=2$. Therefore, $f(2)=2$ or $\alpha$, where $\alpha$ is the unique vertex of $\mathbb{N}_{3}(2)$ in the tournament $T[\mathbb{N}_{3} \cup e \cup e']$.
\epr

\begin{cor}\label{sommets non-critiques}
For all $T\in \mathcal{T}$, we have $V(T)\setminus W_{5}(T)=\sigma(T)$.
\end{cor}
\bpr
Let $T$ be a tournament of the class $\mathcal{T}$ such that $V(T)\setminus W_{5}(T)= \{0,1\}$. By Proposition \ref{3cycle}, we can assume that $T$ is $C_{3}$-critical.
By the same proposition, it suffices to prove that $\{0,1\}\subseteq \sigma(T)$.
By Corollary \ref{abrit 7 sommets}, there is a subset $X$ of $V(T)$ such that $\mathbb{N}_{3} \subset X$ and $T[X]$ is isomorphic to a tournament on 7 vertices of the class $\mathcal{M}$. Suppose for a contradiction that $T$ admits a critical vertex $i \in\{0,1\}$, and let $Y=X\setminus\{i\}$. By Proposition \ref{7 sommets}, $T[Y]$ is indecomposable. As $T$ is $T[Y]$-critical, then $i \not\in$ Ext$(Y)$ by Theorem~\ref{partiellement critique}. This is a contradiction because $T[X]$ is indecomposable.
\epr

Now, we prove that $\mathcal{M} \subseteq \mathcal{T}$. More precisely:
\begin{prop}\label{M dans T}
For all tournament $T$ of the class $\mathcal{M}$, we have $V(T)\setminus W_{5}(T)=\sigma(T)=\{0,1\}$.
\end{prop}
\bpr
Let $T$ be a tournament on $(2n+1)$ vertices of the class $\mathcal{M}$ for some $n\geq3$. By Corollary~\ref{sommets non-critiques}, it suffices to prove that $V(T)\setminus W_{5}(T)=\{0,1\}$. We proceed by induction on $n$. By Proposition \ref{7 sommets}, the statement is satisfied for $n=3$. Let now $n\geq 4$. So, either $T$ is a tournament on 9 vertices of the class $\mathcal{L}\cup \mathcal{L}^{\star}$ or there is $Q\in \mathcal{C}(G^{T}_{\mathbb{N}_{3}})$ such that $\mid\!Q\!\mid\geq 4$. In the first case, for all $e\in E(G^{T}_{\mathbb{N}_{3}})$, $T-e$ is isomorphic to $U_{7}$ or to a tournament on 7 vertices of the class $\mathcal{K}\cup \mathcal{K}^{\star}$. Therefore, if there exists a subset $Z$ of $V(T)$ such that $Z \cap \{0,1\} \neq \emptyset$ and $T[Z]\simeq W_{5}$, then, for all $e\in E(G^{T}_{\mathbb{N}_{3}})$, $e \subset Z$ by Lemma~\ref{G-2n-2}. Thus, $V(T)\setminus \mathbb{N}_{3} \subset Z$, a contradiction. As furthermore, $W_{5}$ embeds into $T$, then $V(T)\setminus W_{5}(T)=\{0,1\}$ by Theorem \ref{HIK}. In the second case, let $Q\in \mathcal{C}(G^{T}_{\mathbb{N}_{3}})$ such that $\mid\!Q\!\mid\geq 4$. Let $\mathcal{X}=\mathcal{H}$, $\mathcal{I}$, $\mathcal{J}$, $\mathcal{K}$ or $\mathcal{L}$. For $T\in \mathcal{X}$, by Lemma \ref{G-2n-2}, there is $e\in E(G^{T}_{\mathbb{N}_{3}}[Q]) $ such that $T-e$ is $C_{3}$-critical. Moreover, $T-e$ is isomorphic to a tournament of the class $\mathcal{X}$ because $\mathcal{C}(G^{T-e}_{\mathbb{N}_{3}})$ is as described in the same class. By induction hypothesis, $W_{5}$ embeds into $T-e$, and thus into $T$. By Theorem \ref{HIK}, it suffices to verify that $\{0,1\}\subseteq V(T)\setminus W_{5}(T)$. So, Suppose that there exists $Z \subset V(T)$ such that $Z \cap \{0,1\} \neq \emptyset$ and $T[Z]\simeq W_{5}$. By induction hypothesis and by Lemma \ref{G-2n-2}, $Q \subset Z$, so that $Z \subset  Q\cup \mathbb{N}_{3}$. This is a contradiction by Theorem \ref{Latka}, because $T[\mathbb{N}_{3} \cup Q]\simeq U_{\mid\!Q\!\mid+3}$ or $T_{\mid\!Q\!\mid+3}$ by Lemma \ref{une composante}.
\epr

We are now ready to construct the tournaments of the class $\mathcal{T}$. We partition these tournaments $T$ according to the following invariant $c(T)$. For $T\in \mathcal{T}$, $c(T)$ is the minimum of $\mid\! \mathcal{C}(G^{T}_{\sigma(T)\cup\{x\}})\!\mid$, the minimum being taken over all the vertices $x$ of $W_{5}(T)$ such that $T[\sigma(T)\cup\{x\}]\simeq C_{3}$. Notice that $c(T)=c(T^{\star})$. As $T$ is $T[\sigma(T)\cup\{x\}]$-critical by Proposition \ref{3cycle}, then $c(T)\leq 4$. Moreover, $c(T)\geq 2 $ by Lemma \ref{une composante}. Proposition \ref{3cycle} leads us to classify the tournaments $T$ of the class $\mathcal{T}$ according to the different values of $c(T)$. We will see that $c(T)=2$ or 3. Theorem \ref{HIK2} results from Propositions \ref{M dans T}, \ref{2 composantes}, \ref{3 composantes} and \ref{2 ou 3}.

\begin{prop}\label{2 composantes}
 Up to isomorphy, the tournaments $T$ of the class $\mathcal{T}$ such that $c(T)=~2$ are those of the class $\mathcal{M}\setminus (\mathcal{L}\cup \mathcal{L}^{\star})$.
\end{prop}
\begin{sloppypar}
\bpr
For all $T\in \mathcal{M}\setminus (\mathcal{L}\cup \mathcal{L}^{\star})$, we have $T\in \mathcal{T}$ by Proposition~\ref{M dans T}, and $c(T)=2$ by Lemma \ref{une composante}.
Now, let $T$ be a tournament on $(2n+1)$ vertices of the class $\mathcal{T}$ such that $c(T)=2$.  By Proposition \ref{3cycle}, we can assume that $T$ is $C_{3}$-critical with $V(T)\setminus W_{5}(T)=\{0,1\}$ and $\mid\! \mathcal{C}(G^{T}_{\mathbb{N}_{3}}) \!\mid = 2$. By Corollary \ref{abrit 7 sommets} and by interchanging $T$ and $T^{\star}$, there is a tournament $R$ on 7 vertices of the class $\mathcal{H} \cup \mathcal{I} \cup \mathcal{J} \cup \mathcal{K}$ such that for all $e \in E(G^{T}_{\mathbb{N}_{3}}[Q])$ and for all $e'\in E(G^{T}_{\mathbb{N}_{3}}[Q'])$, there exists an isomorphism $f$, fixing 0 and 1, from $R$ onto $T[\mathbb{N}_{3} \cup e\cup e']$, where $Q$ and $Q'$ are the two different connected components of $G^{T}_{\mathbb{N}_{3}}$. If $f(2)=2$, then, by Theorem \ref{partiellement critique}, $T$ and $R$ are in the same class $\mathcal{H}$, $\mathcal{I}$, $\mathcal{J}$ or $\mathcal{K}$. Suppose now that $f(2)\neq 2$. By Corollary \ref{abrit 7 sommets}, $R\in \mathcal{I} \cup \mathcal{K}$. If $R\in \mathcal{I}$ (resp. $\mathcal{K}$), then $T[\mathbb{N}_{3} \cup e\cup e']$ is a tournament on 7 vertices of the class $\mathcal{I}'$ (resp. $\mathcal{K}'$) of the $C_{3}$-critical tournaments $Z$ such that $\mathcal{C}(G^{Z}_{\mathbb{N}_{3}})=\{\mathbb{N}_{3}^{-}(0) \cup \mathbb{N}_{3}^{+}(1), \mathbb{N}_{3}^{-}(1) \cup \mathbb{N}_{3}^{-}(2)\}$ (resp. $\mathcal{C}(G^{Z}_{\mathbb{N}_{3}})=\{\mathbb{N}_{3}^{+}(1) \cup \mathbb{N}_{3}^{-},\mathbb{N}_{3}^{-}(1) \cup \mathbb{N}_{3}^{+}(2)\}$). By Theorem \ref{partiellement critique}, $T\in \mathcal{I}'$ (resp. $ \mathcal{K}'$). Moreover, by considering the vertex $\alpha=\min T[\mathbb{N}_{3}^{-}(2)]$ (resp. $\max T[\mathbb{N}_{3}^{+}(2)]$) and by using Corollary \ref{sommets non-critiques}, $T$ is also $T[\{0,1, \alpha \}]$-critical with $\mathcal{C}(G^{T}_{\{0,1, \alpha \}})=\{\{0,1, \alpha \}^{-}(0) \cup \{0,1, \alpha \}^{+}(1), \{0,1, \alpha \}^{+}(0) \cup \{0,1, \alpha \}^{+}(\alpha)\}$ (resp. $\mathcal{C}(G^{T}_{\{0,1, \alpha \}})=\{\{0,1, \alpha \}^{+}(1) \cup \{0,1, \alpha \}^{-}, \{0,1, \alpha \}^{+}(0) \cup \{0,1, \alpha \}^{-}(\alpha)\}$). It follows that $T$ is isomorphic to a tournament of the class $\mathcal{I}$ (resp. $\mathcal{K}$). \epr
\end{sloppypar}

\begin{prop}\label{3 composantes}
Up to isomorphy, the tournaments $T$ of the class $\mathcal{T}$ such that $c(T)=3$ are those of the class $\mathcal{L}\cup \mathcal{L}^{\star}$.
\end{prop}
\begin{sloppypar}
\bpr
Let $T$ be a tournament of the class $\mathcal{L}\cup \mathcal{L}^{\star}$. $T\in \mathcal{T}$ by Proposition \ref{M dans T}. Moreover, $c(T)=3$ by Theorem \ref{partiellement critique}. Indeed, it suffices to observe that for all $x \in \{i \in V(T) \setminus \mathbb{N}_{3}: T[\{0,1,i\}]\simeq C_{3}\} = \mathbb{N}_{3}^{-}(2)$, we have $\max T[ \mathbb{N}_{3}^{+}(1)]\in X^{+}(1)$, $\min T[ \mathbb{N}_{3}^{-}]\in X^{-}$, $\min T[ \mathbb{N}_{3}^{+}]\in X^{+}$, $\max T[ \mathbb{N}_{3}^{-}(0)]\in X^{-}(0)$ and $2\in X^{+}(x)$, where $X = \{0,1,x\}$.

Now, let $T$ be a tournament on $(2n+1)$ vertices of $\mathcal{T}$ such that $c(T)=3$. By Proposition~\ref{3cycle}, we can assume that $T$ is $C_{3}$-critical with $V(T)\setminus W_{5}(T)=\{0,1\}$ and $\mid\! \mathcal{C}(G^{T}_{\mathbb{N}_{3}}) \!\mid = 3$. By Corollary \ref{abrit 7 sommets} and by interchanging $T$ and $T^{\star}$, there is a tournament $R$ on 7 vertices of the class $\mathcal{H} \cup \mathcal{I} \cup \mathcal{J} \cup \mathcal{K}$ such that for all $e \in E(G^{T}_{\mathbb{N}_{3}}[Q])$ and $e'\in E(G^{T}_{\mathbb{N}_{3}}[Q'])$, there exists an isomorphism $f$, which fixes 0 and 1, from $R$ onto $T[\mathbb{N}_{3} \cup e\cup e']$, where $Q \neq Q'\in \mathcal{C}(G^{T}_{\mathbb{N}_{3}})$. Take $e''\in (G^{T}_{\mathbb{N}_{3}}[Q''])$, where $Q''=\mathcal{C}(G^{T}_{\mathbb{N}_{3}})\setminus \{Q,Q'\}$. Suppose, toward a contradiction, that $R\in\mathcal{H} \cup \mathcal{J}$. By Theorem \ref{partiellement critique} and by Corollary~\ref{abrit 7 sommets}, if $R\in\mathcal{H}$ (resp. $R\in\mathcal{J}$), then $\{Q,Q'\}=\{\mathbb{N}_{3}^{+}(0)\cup \mathbb{N}_{3}^{-}, \mathbb{N}_{3}^{-}(1)\cup \mathbb{N}_{3}^{+}\}$ (resp. $\{\mathbb{N}_{3}^{+}(1)\cup \mathbb{N}_{3}^{-}, \mathbb{N}_{3}^{-}(0)\cup \mathbb{N}_{3}^{-}(1)\}$). Therefore, by Lemma \ref{une composante}, $Q''= \{\mathbb{N}_{3}^{+}(1)\cup \mathbb{N}_{3}^{+}(2)\}$, $\{\mathbb{N}_{3}^{-}(0)\cup \mathbb{N}_{3}^{+}(1)\}$ or $ \{\mathbb{N}_{3}^{-}(0)\cup \mathbb{N}_{3}^{-}(2)$\} (resp. $ \{\mathbb{N}_{3}^{+}(0)\cup \mathbb{N}_{3}(2)\}$ or $\{\mathbb{N}_{3}^{+}\cup \mathbb{N}_{3}^{-}(2)\}$). We verify that in each of these cases, either $T[\{0\}\cup e\cup e'']$, $T[\{0\}\cup e'\cup e'']$, $T[\{1\}\cup e\cup e'']$ or $T[\{1\}\cup e'\cup e'']$ is isomorphic to $W_{5}$, a contradiction. Therefore, $R\in \mathcal{I} \cup \mathcal{K} $. By Corollary \ref{abrit 7 sommets}, $f(2)=2$ or $\alpha$, where $\alpha$ is the unique vertex of $\mathbb{N}_{3}(2)$ in $T[\mathbb{N}_{3} \cup e\cup e']$.

Suppose, again by contradiction, that $R\in \mathcal{I}$. We begin by the case where $f(2)=2$. By Theorem \ref{partiellement critique}, we can suppose that $Q=\{\mathbb{N}_{3}^{+}(0)\cup \mathbb{N}_{3}^{+}(2)\}$ and $Q'=\{ \mathbb{N}_{3}^{-}(0)\cup \mathbb{N}_{3}^{+}(1)\}$. By Lemma \ref{une composante}, $Q''= \{\mathbb{N}_{3}^{-}(1)\cup \mathbb{N}_{3}^{+}\}$, $\{\mathbb{N}_{3}^{-}(2)\cup \mathbb{N}_{3}^{+}\}$ or $\{\mathbb{N}_{3}^{-}(1)\cup \mathbb{N}_{3}^{-}(2)\}$. If $Q''= \{\mathbb{N}_{3}^{-}(1)\cup \mathbb{N}_{3}^{+}\}$ (resp. $\{\mathbb{N}_{3}^{-}(2)\cup \mathbb{N}_{3}^{+}\}$), then $T[\{0\}\cup e\cup e'']\simeq W_{5}$ (resp. $T[\{1\}\cup e'\cup e'']\simeq W_{5}$), a contradiction. If $Q''=\{\mathbb{N}_{3}^{-}(1)\cup \mathbb{N}_{3}^{-}(2)\}$, then, by taking $X = \{0,1,x\}$, where $x = \min T[\mathbb{N}_{3}^{-}(2)]$, we obtain a contradiction because, by Corollary \ref{sommets non-critiques}, $T$ is $T[X]$-critical with $\mid\! \mathcal{C}(G^{T}_{X}) \!\mid = 2$. Indeed, $\mathcal{C}(G^{T}_{X}) = \{X^{-}(0) \cup X^{+}(1), X^{+}(0) \cup X^{+}(x) \}$, with $X^{-}(0) = \mathbb{N}_{3}^{-}(0)$, $X^{+}(1) = \mathbb{N}_{3}^{+}(1)$, $X^{+}(0) = \mathbb{N}_{3}^{+}(0) \cup \mathbb{N}_{3}^{-}(1)$ and $X^{+}(x) = \mathbb{N}_{3}^{+}(2) \cup \{2\} \cup (\mathbb{N}_{3}^{-}(2) \setminus \{x\})$. Now, if $f(2)=\alpha$, then we obtain again a contradiction. Indeed, by replacing $T$ by $T^{\star}$ and by interchanging the vertices 0 and 1, $\{Q,Q'\}=\{\mathbb{N}_{3}^{+}(0)\cup \mathbb{N}_{3}^{+}(2), \mathbb{N}_{3}^{-}(0)\cup \mathbb{N}_{3}^{+}(1)\}$ as in the case where $f(2)=2$.

At present, $R\in \mathcal{K}$. We begin by the case where $f(2)=2$. By Theorem \ref{partiellement critique}, we can suppose that  $Q=\{\mathbb{N}_{3}^{-}\cup \mathbb{N}_{3}^{+}(1)\}$ and $Q'=\{\mathbb{N}_{3}^{+}(0)\cup \mathbb{N}_{3}^{-}(2)\}$. By Lemma~\ref{une composante}, $Q''=\{ \mathbb{N}_{3}^{-}(0)\cup \mathbb{N}_{3}^{+}\}$, $\{\mathbb{N}_{3}^{-}(1)\cup \mathbb{N}_{3}^{+}\}$, $\{\mathbb{N}_{3}^{-}(0)\cup \mathbb{N}_{3}^{-}(1)\}$ or $\{\mathbb{N}_{3}^{-}(1)\cup \mathbb{N}_{3}^{+}(2)\}$. If $Q''=\{ \mathbb{N}_{3}^{-}(1)\cup \mathbb{N}_{3}^{+}\}$ (resp. $\{\mathbb{N}_{3}^{-}(1)\cup \mathbb{N}_{3}^{-}(0)\}$), then $T[\{0\}\cup e\cup e'']\simeq W_{5}$ (resp. $T[\{1\}\cup e'\cup e'']\simeq W_{5}$), a contradiction. If $Q''=\{\mathbb{N}_{3}^{-}(1)\cup \mathbb{N}_{3}^{+}(2)\}$, then, by taking $X = \{0,1,x\}$, where $x = \max T[\mathbb{N}_{3}^{+}(2)]$, we have a contradiction because, by Corollary \ref{sommets non-critiques}, $T$ is $T[X]$-critical with $\mid\! \mathcal{C}(G^{T}_{X}) \!\mid = 2$. Indeed, $\mathcal{C}(G^{T}_{X}) = \{X^{-} \cup X^{+}(1), X^{+}(0) \cup X^{-}(x) \}$, with $X^{-} = \mathbb{N}_{3}^{-}$, $X^{+}(1) = \mathbb{N}_{3}^{+}(1)$, $X^{+}(0) = \mathbb{N}_{3}^{-}(1) \cup \mathbb{N}_{3}^{+}(0)$ and $X^{-}(x) = \mathbb{N}_{3}^{-}(2) \cup \{2\} \cup (\mathbb{N}_{3}^{+}(2) \setminus \{x\})$. If $Q''=\{\mathbb{N}_{3}^{-}(0)\cup \mathbb{N}_{3}^{+}\}$, then $T\in \mathcal{L}$. Now, suppose that $f(2)=\alpha$.  By Theorem \ref{partiellement critique}, we can suppose that $Q=\{\mathbb{N}_{2}^{+}(1) \cup \mathbb{N}_{2}^{-}\}$ and $Q'=\{\mathbb{N}_{2}^{-}(1) \cup \mathbb{N}_{3}^{+}(2)\}$. By Lemma \ref{une composante}, $Q''= \{\mathbb{N}_{3}^{-}(2)\cup \mathbb{N}_{3}^{+}\}$, $ \{\mathbb{N}_{3}^{-}(2)\cup \mathbb{N}_{3}(0)\}$ or $\{\mathbb{N}_{3}^{-}(0)\cup \mathbb{N}_{3}^{+}\}$. If $Q''= \{\mathbb{N}_{3}^{-}(2)\cup \mathbb{N}_{3}^{+}\}$ or $\{\mathbb{N}_{3}^{-}(0)\cup \mathbb{N}_{3}^{-}(2)\}$, then $T[\{0\}\cup e\cup e'']\simeq W_{5}$, a contradiction. If $Q''=\{\mathbb{N}_{3}^{+}(0)\cup \mathbb{N}_{3}^{-}(2)\}$, then we obtain the same configuration giving $\mid\! \mathcal{C}(G^{T}_{X}) \!\mid = 2$ in the case where $f(2)=2$.  If $Q''=\{\mathbb{N}_{3}^{-}(0)\cup \mathbb{N}_{3}^{+}\}$, then $T$ is isomorphic to a tournament of the class $\mathcal{L}^{\star}$. \epr

\end{sloppypar}

\begin{prop}\label{2 ou 3}
For any tournament $T$ of the class $\mathcal{T}$, we have $c(T)=2$ or 3.
\end{prop}
\bpr
Let $T$ be a tournament on $(2n+1)$ vertices of the class $\mathcal{T}$ for some $n\geq 3$. We proceed by induction on $n$. By Propositions \ref{2 composantes} and \ref{3 composantes}, the statement is satisfied for $n=3$ and for $n=4$. Let $n\geq 5$. By Proposition \ref{3cycle}, we can assume that $T$ is $C_{3}$-critical with $V(T)\setminus W_{5}(T)=\{0,1\}$. By Theorem \ref{partiellement critique} and Lemma~\ref{une composante}, $2 \leq c(T) \leq 4$. So, we only consider the case where $\mid\!\mathcal{C}(G^{T}_{\mathbb{N}_{3}})\!\mid=4$. By Corollary \ref{abrit 7 sommets}, there exist $Q\neq Q'\in \mathcal{C}(G^{T}_{\mathbb{N}_{3}})$ and a tournament $R$ on 7 vertices of the class $\mathcal{M}$, such that for all $e \in E(G^{T}_{\mathbb{N}_{3}}[Q])$ and for all $e'\in E(G^{T}_{\mathbb{N}_{3}}[Q'])$, $T[\mathbb{N}_{3} \cup e\cup e']\simeq R$.
By Lemma \ref{G-2n-2}, there exists $e''\in E(G^{T}_{\mathbb{N}_{3}}[Q''])$, where $Q''\in \mathcal{C}(G^{T}_{\mathbb{N}_{3}})\setminus \{Q,Q'\}$, such that $T-e''$ is $C_{3}$-critical. As $W_{5}$ embeds into $T-e''$, then $ V(T-e'')\setminus W_{5}(T-e'')=\{0,1\}$ by Theorem \ref{HIK}. Therefore, $T-e''\in \mathcal{T}$. By induction hypothesis, $c(T-e'')=2$ or 3. By Theorem \ref{partiellement critique}, if $c(T-e'')=2$, then $c(T)=2$ or 3. So, suppose that $c(T-e'')=3$. By Proposition~\ref{3 composantes} and by interchanging $T$ and $T^{\star}$, we can assume that $T-e''\in\mathcal{L}$. By Theorem \ref{partiellement critique} and by taking $e''=\{x,x'\}$, we can assume that $x\in \mathbb{N}_{3}^{-}(1)$ and $x'\in \mathbb{N}_{3}^{+}(2)$. Thus, for $X=\{0,1,x'\}$, we have $T[X] \simeq C_{3}$ and $X^{+}(x')=\emptyset$. It follows by Theorem \ref{partiellement critique} that $c(T) < 4$.
\epr

\end{document}